\documentclass{elsart}
\usepackage{amssymb}
\usepackage{graphicx}

\begin{document}

\begin{frontmatter}

\title{Global bifurcations of limit cycles
in~a~Holling-type~dynamical~system\thanksref{label1}}
\thanks[label1]{This work was partially supported by grants from
the Simons Foundation (International Mathematical Union) and
the Netherlands Organization for Scientific Research (NWO).
The author is also very grateful to the Department of Ma\-the\-ma\-tics
and Statistics of the Missouri University of Science and Technology (USA)
for hospitality and support during his stay in March\,--\,April 2015.}

\author{Valery A. Gaiko}

\ead{valery.gaiko@gmail.com}

\address{National~Academy~of~Sciences~of~Belarus,
United~Institute~of~Informatics~Problems,
Minsk~220012,~Belarus}

\begin{abstract}
In this paper, we complete the global qualitative analysis of a quartic family
of planar vector fields corresponding to a rational Holling-type dynamical system
which models the dynamics of the populations of predators and their prey in a given
ecological or biomedical system. In particular, studying global bifurcations, we prove
that such a system can have at most two limit cycles surrounding one singular point.
\end{abstract}

\end{frontmatter}

\section{Introduction}

In this paper, we consider a quartic family of planar vector fields corresponding
to a rational Holling-type dynamical system which models the dynamics of the populations
of predators and their prey in a given ecological or biomedical system and which is
a variation on the classical Lotka--Volterra system. For the latter system the change
of the prey density per unit of time per predator called the response function is
proportional to the prey density. This means that there is no saturation of the predator
when the amount of available prey is large. However, it is more realistic to consider
a nonlinear and bounded response function, and in fact different response functions
have been used in the literature to model the predator response; see
\cite{Bazykin}--\cite{bg}, \cite{hol}--\cite{lll}, \cite{zcw}.
    \par
For instance, in \cite{zcw}, the following predator--prey model has been studied:
    $$
    \begin{array}{l}
\dot{x}=x(a-\lambda x)-yp(x) \qquad \mbox{(prey)},
    \\[2mm]
\dot{y}=-\delta y+yq(x) \ \ \ \qquad \mbox{(predator)}.
    \\[2mm]
    \end{array}
    \eqno(1.1)
    $$
The variables $x>0$ and $y>0$ denote the density of the prey and predator populations
respectively, while $p(x)$ is a non-monotonic response function given by
    $$
p(x)=\frac{mx}{\alpha x^{2}+\beta x+1},
    \\[2mm]
    \eqno(1.2)
    $$
where $\alpha,$ $m$ are positive and where $\beta>-2\sqrt{\alpha}.$ Observe that in the
absence of predators, the number of prey increases according to a logistic growth law.
The coefficient a represents the intrinsic growth rate of the prey, while $\lambda>0$
is the rate of competition or resource limitation of prey. The natural death rate of
the predator is given by $\delta>0.$ In Gause's model the function $q(x)$ is given by
$q(x)=cp(x),$ where $c>0$ is the rate of conversion between prey and predator \cite{zcw}.
    \par
In \cite{bnrs}, \cite{bg}, the following family has been investigated:
    $$
    \begin{array}{l}
\displaystyle\dot{x}=x\left(1-\lambda x-\frac{y}{\alpha x^{2}+\beta x+1}\right),
    \\[4mm]
\displaystyle\dot{y}=-y\left(\delta+\mu y-\frac{x}{\alpha x^{2}+\beta x+1}\right),
    \end{array}
    \eqno(1.3)
    $$
where $\alpha\geq0,$ $\beta>-2\sqrt{\alpha},$ $\delta>0,$ $\lambda>0,$ and $\mu\geq0$
are parameters. Note that (1.3) is obtained from (1.1) by adding the term $-\mu y^{2}$
to the second equation and after scaling $x$ and $y,$ as well as the parameters and
the time $t.$ In this way, it has been taken into account competition between predators
for resources other than prey. The non-negative coefficient $\mu$ is the rate of
competition amongst predators. Systems (1.1)--(1.3) represent predator--prey models
with a generalized Holling response functions of type~IV.
    \par
In \cite{lcr}, it has been considered the following generalized Gause predator--prey system
    $$
    \begin{array}{l}
\dot{x}=rx(1-x/k)-yp(x),
    \\[2mm]
\dot{y}=y(-d+cp(x))
    \\[2mm]
    \end{array}
    \eqno(1.4)
    $$
with a generalized Holling response function of type~III:
    $$
p(x)=\frac{mx^{2}}{ax^{2}+bx+1}.
    \eqno(1.5)
    $$
This system, where $x>0$ and $y>0,$ has seven parameters: the parameters $a,$ $c,$ $d,$ $k,$
$m,$ $r$ are positive and the parameter $b$ can be negative or non-negative. The parameters
$a,$ $b,$ and $m$  fitting parameters of response function. The parameter $d$ is the death rate
of the predator while $c$ is the efficiency of the predator to convert prey into predators.
The prey follows a logistic growth with a rate $r$ in the absence of predator. The environment
has a prey capacity determined by $k.$
    \par
The case $b\geq0$ has been studied earlier; see the references in \cite{lcr}. The case $b<0$
is more interesting: it provides a model for a functional response with limited group defence.
In opposition to the generalized Holling function of type~IV studied in \cite{bnrs}, \cite{bg}, \cite{zcw}, where the response function tends to zero as the prey population tends to infinity,
the generalized function of type~III tends to a non-zero value as the prey population tends
to infinity. The functional response of type III with $b<0$ has a maximum at some point;
see \cite{lcr}. When studying the case $b<0,$ one can find also a Bogdanov--Takens bifurcation of
codimension~3 which is an organizing center for the bifurcation diagram of system (1.10)--(1.5)
\cite{lcr}.
    \par
After scaling $x$ and $y,$ as well as the parameters and the time $t,$ this system can be
reduced to a system with only four parameters $(\alpha,$ $\beta,$ $\delta,$ $\rho)$ \cite{lcr}:
    $$
    \begin{array}{l}
\dot{x}=\rho x(1-x)-yp(x),
    \\[2mm]
\dot{y}=y(-\delta+p(x)),
    \\[-2mm]
    \end{array}
    \eqno(1.6)
    $$
where
    $$
p(x)=\frac{x^{2}}{\alpha x^{2}+\beta x+1}.
    \eqno(1.7)
    $$
    \par
Note that many studies of discrete-time predator--prey models have been done;
see, e.\,g., \cite{lll} and the references therein. In particular, the following
discrete predator--prey system with a generalized Holling response function
of type~III has been investigated in \cite{lll}:
    $$
    \begin{array}{l}
\displaystyle x_{k+1}=x_{k}+\theta x_{k}\left(\rho(1-x_{k})(x_{k}-\lambda)-\frac{x_{k}y_{k}}
{\alpha x_{k}^{2}+\beta x_{k}+1}\right),
    \\[4mm]
\displaystyle y_{k+1}=y_{k}+\theta y_{k}\left(-\delta+\frac{\gamma x_{k}^{2}}
{\alpha x_{k}^{2}+\beta x_{k}+1}\right),
    \end{array}
    \eqno(1.8)
    $$
where $\theta>0$ is the step size in the forward Euler scheme.
    \par
In this paper, we study the system
    $$
    \begin{array}{l}
\displaystyle\dot{x}=x\left(1-\lambda x-\frac{xy}{\alpha x^{2}+\beta x+1}\right)
\qquad \qquad \mbox{(prey)},
    \\[4mm]
\displaystyle\dot{y}=-y\left(\delta+\mu y-\frac{x^{2}}{\alpha x^{2}+\beta x+1}\right)
\qquad \mbox{(predator)},
    \end{array}
    \eqno(1.9)
    $$
where $x>0$ and $y>0;$ $\alpha\geq0,$ $-\infty<\beta<+\infty,$ $\delta>0,$ $\lambda>0,$
and $\mu\geq0$ are parameters.
    \par
Rational system (1.9) can be written in the form of a quartic dynamical system
    $$
    \begin{array}{l}
\dot{x}=~~x((1-\lambda x)(\alpha x^{2}+\beta x+1)-xy)\equiv P,
    \\[2mm]
\dot{y}=-y((\delta +\mu y)(\alpha x^{2}+\beta x+1)-x^{2})\equiv Q.
    \end{array}
    \eqno(1.10)
    $$
Together with (1.10), we will also consider an auxiliary system
(see \cite{BL}, \cite{Perko})
    $$
\dot{x}=P-\gamma Q, \qquad \dot{y}=Q+\gamma P,
    \eqno(1.11)
    $$
applying to these systems new bifurcation methods and geometric approaches developed
in \cite{bg}--\cite{gai10} and completing the qualitative analysis of (1.9).

\section{Basic facts on singular points and limit cycles}

The study of singular points of system (1.9) will use two index theorems by H.\,Poincar\'{e};
see \cite{BL}. But first let us define the singular point and its Poincar\'{e} index~\cite{BL}.
    \medskip
    \par
    \textbf{Definition 2.1.}
A singular point of the dynamical system
    $$
    \dot{x}=P(x,y), \quad \dot{y}=Q(x,y),
    \eqno(2.1)
    $$
where $P(x,y)$ and $Q(x,y)$ are continuous functions (for example, polynomials),
is a point at which the right-hand sides of (2.1) simultaneously vanish.
    \medskip
    \par
    \textbf{Definition 2.2.}
Let $S$ be a simple closed curve in the phase plane not passing through a singular
point of system (2.1) and $M$ be some point on $S.$ If the point $M$ goes around
the curve $S$ in positive direction (counterclockwise) one time,
then the vector coinciding with the direction of a tangent to
the trajectory passing through the point $M$ is rotated through
the angle $2\pi j$ $(j=0,\pm1,\pm2,\ldots).$ The integer $j$
is called the \emph{Poincar\'{e} index} of the closed curve $S$
relative to the vector field of system~(2.1) and has the expression
    $$
    j=\frac{1}{2\pi}\oint_S\frac{P~dQ-Q~dP}{P^2+Q^2}.\\[-2mm]
    $$
    \par
According to this definition, the index of a node or a focus, or a
center is equal to $+1$ and the index of a saddle is $-1.$
    \par
    \medskip
    \textbf{Theorem 2.1 (First Poincar\'{e} Index Theorem).}
    \emph{If $N,$ $N_f,$ $N_c,$ and $C$ are respectively the number
of nodes, foci, centers, and saddles in a finite part of the phase
plane and $N'$ and $C'$ are the number of nodes and saddles at
infinity, then it is valid the formula}
    $$
    N+N_f+N_c+N'=C+C'+1.
    $$
    \par
    \textbf{Theorem 2.2 (Second Poincar\'{e} Index Theorem).}
    \emph{If all singular points are simple, then along an isocline
without multiple points lying in a Poincar\'{e} hemisphere which is
obtained by a stereographic projection of the phase plane, the
singular points are distributed so that a saddle is followed by
a node or a focus, or a center and vice versa. If two points are
separated by the equator of the Poincar\'{e} sphere, then a saddle
will be followed by a saddle again and a node or a focus, or
a center will be followed by a node or a focus, or a center.}
    \medskip
    \par
Consider polynomial system (2.1) in the vector form
    $$
    \mbox{\boldmath$\dot{x}$}=\mbox{\boldmath$f$}
    (\mbox{\boldmath$x$},\mbox{\boldmath$\mu$)},
    \eqno(2.2)
    $$
where $\mbox{\boldmath$x$}\in\textbf{R}^2;$ \
$\mbox{\boldmath$\mu$}\in\textbf{R}^n;$ \
$\mbox{\boldmath$f$}\in\textbf{R}^2$ \ $(\,\mbox{\boldmath$f$}$
is a polynomial vector function).
    \par
Let us recall some basic facts concerning limit cycles of (2.2).
But first of all, let us state two fundamental theorems from theory
of ana\-ly\-tic functions~\cite{Gaiko}.
    \par
    \medskip
    \textbf{Theorem 2.3 (Weierstrass Preparation Theorem).}
    \emph{Let $F(w,z)$ be an analytic in the neighborhood of the point
$(0,0)$ function satisfying the following conditions}
    $$
    F(0,0)=0, \: \frac{\partial F(0,0)}{\partial w}=0, \:
    \ldots, \: \frac{\partial^{k-1}F(0,0)}{\partial^{k-1}w}=0; \quad
    \frac{\partial^{k}F(0,0)}{\partial^{k}w}\neq0.
    $$
    \par
    \emph{Then in some neighborhood $|w|<\varepsilon,$ $|z|<\delta$ of
the points $(0,0)$ the function $F(w,z)$ can be represented as}
    \vspace{-2mm}
    $$
    F(w,z)=(w^{k}+A_{1}(z)w^{k-1}+\ldots+A_{k-1}(z)w+A_{k}(z))\Phi(w,z),
    $$
\emph{where $\Phi(w,z)$ is an analytic function not equal to zero
in the chosen neighborhood and $A_{1}(z),\ldots,A_{k}(z)$ are
analytic functions for $|z|<\delta.$}
    \par
    \medskip
From this theorem it follows that the equation $F(w,z)=0$ in a sufficiently
small neighborhood of the point $(0,0)$ is equivalent to the equation
    $$
    w^{k}+A_{1}(z)w^{k-1}+\ldots+A_{k-1}(z)w+A_{k}(z)=0,
    $$
which left-hand side is a polynomial with respect to~$w.$ Thus,
the Weierstrass preparation theorem reduces the local study of the
general case of implicit function $w(z),$ defined by the
equation $F(w,z)=0,$ to the case of implicit function,
defined by the algebraic equation with respect to~$w.$
    \par
    \medskip
    \textbf{Theorem 2.4 (Implicit Function Theorem).}
    \emph{Let $F(w,z)$ be an analytic function in the neighborhood of
the point $(0,0)$ and $F(0,0)\!=\!0,$ $F'_{w}(0,0)\!\neq\!0.$}
    \par
    \emph{Then there exist $\delta>0$ and $\varepsilon>0$ such that
for any~$z$ satisfying the condition $|z|<\delta$ the equation
$F(w,z)=0$ has the only solution $w=f(z)$ sa\-tisfying the
condition $|f(z)|<\varepsilon.$ The func\-tion $f(z)$ is
expanded into the series on positive integer powers of~$z$ which
converges for $|z|<\delta,$ i.\,e., it is a single-valued analytic
function of~$z$ which vanishes at $z=0.$}
    \medskip
    \par
Assume that system (2.2) has a limit cycle
    $$
    L_0:\mbox{\boldmath$x$}=\mbox{\boldmath$\varphi$}_0(t)
    $$
of minimal period $T_0$ at some parameter value
$\mbox{\boldmath$\mu$}=\mbox{\boldmath$\mu$}_0\in\textbf{R}^n$
(Fig.~1).
    \par
    \begin{figure}[htb]
\begin{center}
    \includegraphics[width=138.5mm]{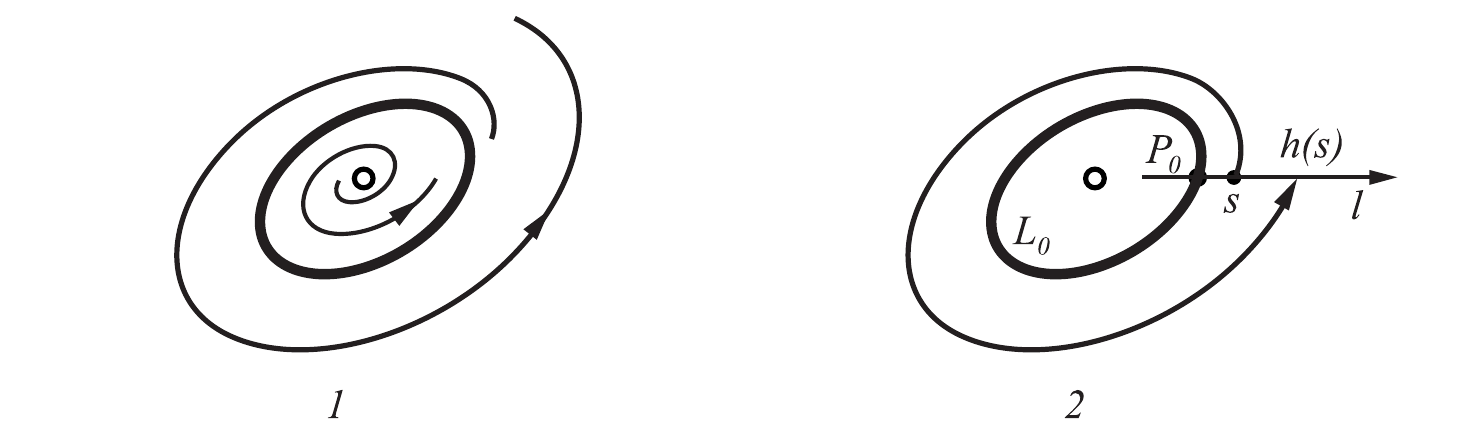}
    \vspace{-4mm}
    \par
    {\small FIG.~1. The Poincar\'{e} return map
    in the neighborhood of a multiple limit cycle.}
\end{center}
    \end{figure}
    \medskip
    \par
Let~$l$ be the straight line normal to $L_0$ at the
point $\mbox{\boldmath$p$}_0=\mbox{\boldmath$\varphi$}_0(0)$ and
$s$ be the coordinate along $l$ with $s$ positive exterior of
$L_0.$ It then follows from the implicit function theorem that
there is a $\delta>0$ such that the Poincar\'{e} map
$h(s,\mbox{\boldmath$\mu$})$ is defined and analytic for $|s|<\delta$
and $\|\mbox{\boldmath$\mu$}-\mbox{\boldmath$\mu$}_0\|<\delta.$
Besides, the displacement function for system (2.2) along
the normal line~$l$ to~$L_0$ is defined as the function
    $$
    d(s,\mbox{\boldmath$\mu$})=h(s,\mbox{\boldmath$\mu$})-s.
    $$
    \par
In terms of the displacement function, a multiple limit cycle
can be defined as follows \cite{Gaiko}.
    \par
    \medskip
    \textbf{Definition 2.3.}
A limit cycle $L_0$ of (2.2) is a \emph{multiple limit cycle} iff
$d(0,\mbox{\boldmath$\mu$}_0)\!=\!d_r(0,\mbox{\boldmath$\mu$}_0)\!=\!0$
and it is a \emph{simple limit cycle} (or hyperbolic limit cycle)
if it is not a multiple limit cycle; furthermore, $L_0$ is a limit
cycle of multiplicity~$m$ iff
    $$
    d(0,\mbox{\boldmath$\mu$}_0)=d_r(0,\mbox{\boldmath$\mu$}_0)=\ldots
    =d_r^{(m-1)}(0,\mbox{\boldmath$\mu$}_0)=0, \quad
    d_r^{(m)}(0,\mbox{\boldmath$\mu$}_0)\neq 0.
    $$
    \par
Note that the multiplicity of $L_0$ is independent of the point
$\mbox{\boldmath$p$}_0\in L_0$ through which we take the normal
line~$l.$
    \par
Let us write down also the following formulas which have already
become classical ones and determine the derivatives of the
displacement function in terms of integrals of the vector
field~$\mbox{\boldmath$f$}$ along the periodic orbit
$\mbox{\boldmath$\varphi$}_0(t)$~\cite{Gaiko}:
    $$
    d_s(0,\mbox{\boldmath$\mu$}_0)\;=\;\displaystyle\exp\int_{0}^{T_0}\!
    \mbox{\boldmath$\nabla$}\cdot\mbox{\boldmath$f$}
    (\mbox{\boldmath$\varphi$}_0(t),\mbox{\boldmath$\mu$}_0)\:\textrm{d}t-1
    \eqno(3.2)
    \vspace{2mm}
    $$
and
    $$
    d_{\mu_j}(0,\mbox{\boldmath$\mu$}_0)=\\
    $$
    $$
    \frac{-\omega\,_0}
    {\|\mbox{\boldmath$f$}(\mbox{\boldmath$\varphi$}_0(0),
    \mbox{\boldmath$\mu$}_0)\|}\;
    \displaystyle\int_{0}^{T_0}\!\!
    \exp\left(-\!\int_{0}^{t}\!\mbox{\boldmath$\nabla$}\cdot
    \mbox{\boldmath$f$}(\mbox{\boldmath$\varphi$}_0(\tau),
    \mbox{\boldmath$\mu$}_0)\,\textrm{d}\tau\right)
    \mbox{\boldmath$f$}\wedge\mbox{\boldmath$f$}_{\mu_j}
    (\mbox{\boldmath$\varphi$}_0(t),\mbox{\boldmath$\mu$}_0)\:\textrm{d}t
    \\[2mm]
    \eqno(3.3)
    $$
for $j=1,\ldots,n,$ where $\omega_0=\pm1$ according to whether $L_0$
is positively or negatively oriented, respectively, and where the wedge
product of two vectors $\mbox{\boldmath$x$}=(x_1,x_2)$ and
$\mbox{\boldmath$y$}=(y_1,y_2)$ in $\textbf{R}^2$ is defined as
    $$
    \mbox{\boldmath$x$}\wedge\mbox{\boldmath$y$}=x_1\,y_2-x_2\,y_1.
    \vspace{-1mm}
    $$
    \par
Similar formulas for $d_{ss}(0,\mbox{\boldmath$\mu$}_0)$ and
$d_{s{\mu_j}}(0,\mbox{\boldmath$\mu$}_0)$ can be derived in terms
of integrals of the vector field $\mbox{\boldmath$f$}$ and its first
and second partial derivatives along $\mbox{\boldmath$\varphi$}_0(t).$
The hypotheses of theorems in the next section will be stated in terms
of conditions on the displacement function $d(s,\mbox{\boldmath$\mu$})$
and its partial derivatives at $(0,\mbox{\boldmath$\mu$}_0)$~\cite{Gaiko}.

\section{Local bifurcation surfaces and the global termination principle
for multiple limit cycles}

In this section, we restate first Perko's theorems on the local existence
of $(n\!-\!m\!+\!1)$-dimensional surfaces, $C_m,$ of multiplicity-$m$
limit cycles for the polynomial system (2.2) with
$\mbox{\boldmath$\mu$}\in\textbf{R}^n$ and $n\geq m\geq2.$
These results describe the topological structure of the codimension
$(m\!-\!1)$ bifurcation surfaces $C_m.$ For $m=2,3,4,$
$C_2,$ $C_3,$ and $C_4$ are the familiar fold, cusp, and
swallow-tail bi\-fur\-ca\-tion surfaces; for $m\geq5,$ the
topological structure of the surfaces~$C_m$ is more complex. For
instance, $C_5$ and $C_6$ are the butterfly and wigwam bifurcation
surfaces, respectively~\cite{Perko}. Since the proofs of the theorems
in this section, describing the universal unfolding near a multiple
limit cycles of~(2.2), pa\-ral\-lel the classical proofs of Catastrophe
Theory, we will only state the theorems (see \cite{Perko} for more detail).
    \medskip
    \par
    \textbf{Definition 3.1.}
An $(n\!-\!1)$-dimensional analytic surface
$C_{2}\subset\textbf{R}^n$ is an \emph{$(n\!-\!1)$-dimensional
fold bifurcation surface of multiplicity-two limit cycles of
(2.2) through a point $\mbox{\boldmath$\mu$}_0\in\textbf{R}^n,$}
if for all $\varepsilon>0$ there exists a $\delta>0$ such that
for each $\mbox{\boldmath$\mu$}\in C_{2}$ with
$\|\mbox{\boldmath$\mu$}-\mbox{\boldmath$\mu$}_{0}\|<\delta,$
the system (2.2) has a unique multiplicity-two limit cycle
$L_{\mbox{\boldmath$\mu$}}$ in an $\varepsilon$-neighborhood
of $L_{0}$ and the system (2.2) undergoes a fold bifurcation
at $L_{\mbox{\boldmath$\mu$}};$ i.\,e., for
$\|\mbox{\boldmath$\mu$}-\mbox{\boldmath$\mu$}_{0}\|<\delta,$
$L_{\mbox{\boldmath$\mu$}}$ splits into a simple stable and a
simple unstable limit cycles in an $\varepsilon$-neighborhood of
$L_{0}$ for $\mbox{\boldmath$\mu$}$ on one side of $C_{2}$ and
$L_{\mbox{\boldmath$\mu$}}$ vanishes for $\mbox{\boldmath$\mu$}$
on the other side of $C_{2}.$ \ Cf. Fig.~2.
    \par
\begin{figure}[htb]
\begin{center}
\includegraphics[width=105mm]{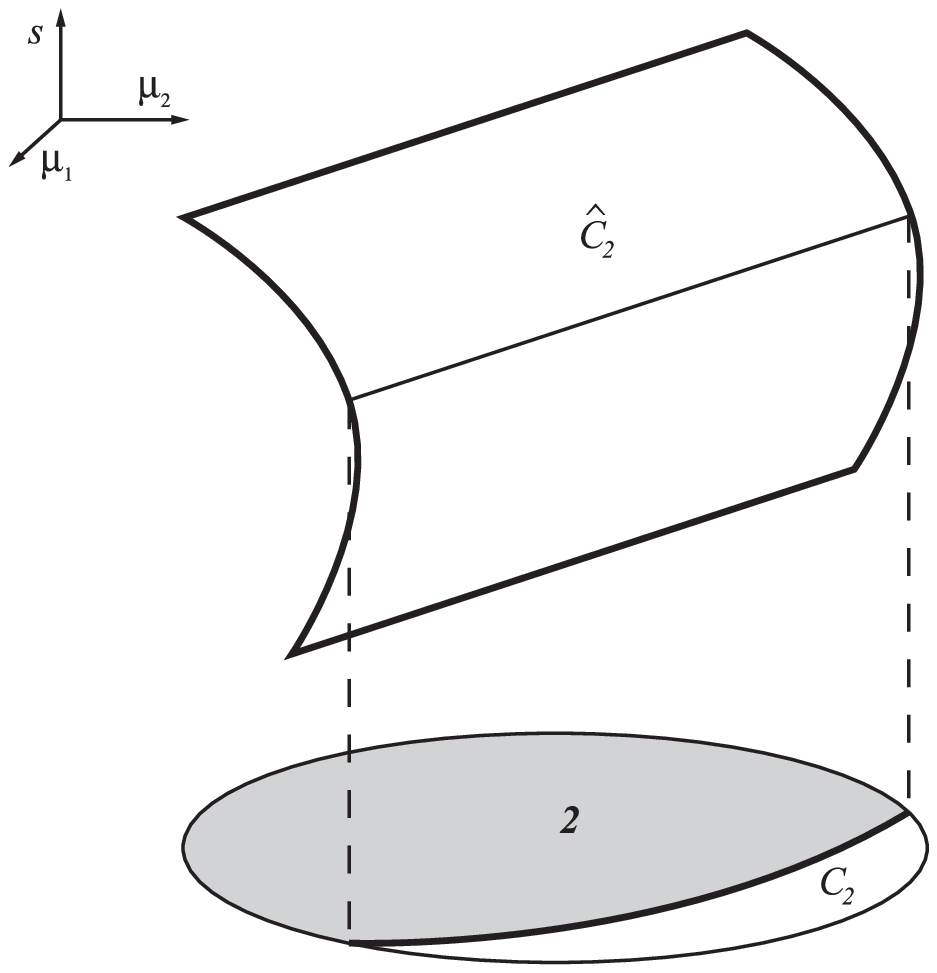}
\vspace{2mm}
    \par
{\small FIG.~2. The fold bifurcation surface.}
\end{center}
\end{figure}
    \par
    \medskip
    \textbf{Theorem 3.1.}
    \emph{Suppose that $n\geq 2,$ that for
    $\mbox{\boldmath$\mu$}=\mbox{\boldmath$\mu$}_0\in\textbf{R}^n$
system~(2.2) has a multiplicity-two limit cycle $L_0,$ and that
$d_{\mu_1}(0,\mbox{\boldmath$\mu$}_0)\neq~0.$}
\emph{Then given $\varepsilon>0,$ there is a $\delta>0$
and a unique function $g(\mu_2,\ldots,\mu_n)$ with
$g(\mu_2^{(0)},\ldots,\mu_n^{(0)})=\mu_1^{(0)},$ defined
and analytic for $\vert\mu_2-\mu_2^{(0)}\vert<\delta,$
$\ldots,\vert\mu_n-\mu_n^{(0)}\vert<\delta,$ such that for
$\vert\mu_2-\mu_2^{(0)}\vert<\delta,\ldots,
\vert\mu_n-\mu_n^{(0)}\vert<\delta,$}
\vspace{-1mm}
    $$
    C_2: \quad \mu_1=g(\mu_2,\ldots,\mu_n)\\[-1mm]
    $$
\emph{is an $(n-1)$-dimensional, analytic fold bifurcation surface
of multiplicity-two limit cycles of~(2.2) through the point
$\mbox{\boldmath$\mu$}_0.$}
    \par
    \medskip
    \textbf{Definition 3.2.}
An analytic surface $C_{3}\subset\textbf{R}^n$ is an
\emph{$(n\!-\!2)$-dimen\-si\-onal cusp bifurcation surface of
multiplicity-three limit cycles of (2.2) through a point
$\mbox{\boldmath$\mu$}_0\in\textbf{R}^n,$} if for all
$\varepsilon>0$ there exists a $\delta>0$ such that
for each $\mbox{\boldmath$\mu$}\in C_{3}$ with
$\|\mbox{\boldmath$\mu$}-\mbox{\boldmath$\mu$}_{0}\|<\delta,$
system (2.2) has a unique multiplicity-three limit cycle
$L_{\mbox{\boldmath$\mu$}}$ in an $\varepsilon$-neighborhood
of $L_{0}$ and the system (2.2) undergoes a cusp bifurcation
at $L_{\mbox{\boldmath$\mu$}};$ i.\,e., $C_{3}$ is the intersection
of two $(n\!-\!1)$-dimensional fold bifurcation surfaces of
multiplicity-two limit cycles of (2.2), $C_{2}^{\pm},$ which
intersect in a cusp along $C_{3};$ for
$\|\mbox{\boldmath$\mu$}-\mbox{\boldmath$\mu$}_{0}\|<\delta$
and for $\mbox{\boldmath$\mu$}$ in the cuspidal region between
$C_{2}^{+}$ and $C_{2}^{-}$ (shaded in Fig.~3), system (2.2)
has three simple limit cycles in an $\varepsilon$-neighborhood of $L_{0};$
for $\|\mbox{\boldmath$\mu$}-\mbox{\boldmath$\mu$}_{0}\|<\delta$ and
$\mbox{\boldmath$\mu$}$ outside the cuspidal region, system (2.2) has
one simple limit cycle in an $\varepsilon$-neighborhood of~$L_{0}.$
Cf.~Fig.~3.
    \par
\begin{figure}[htb]
\begin{center}
\includegraphics[width=105mm]{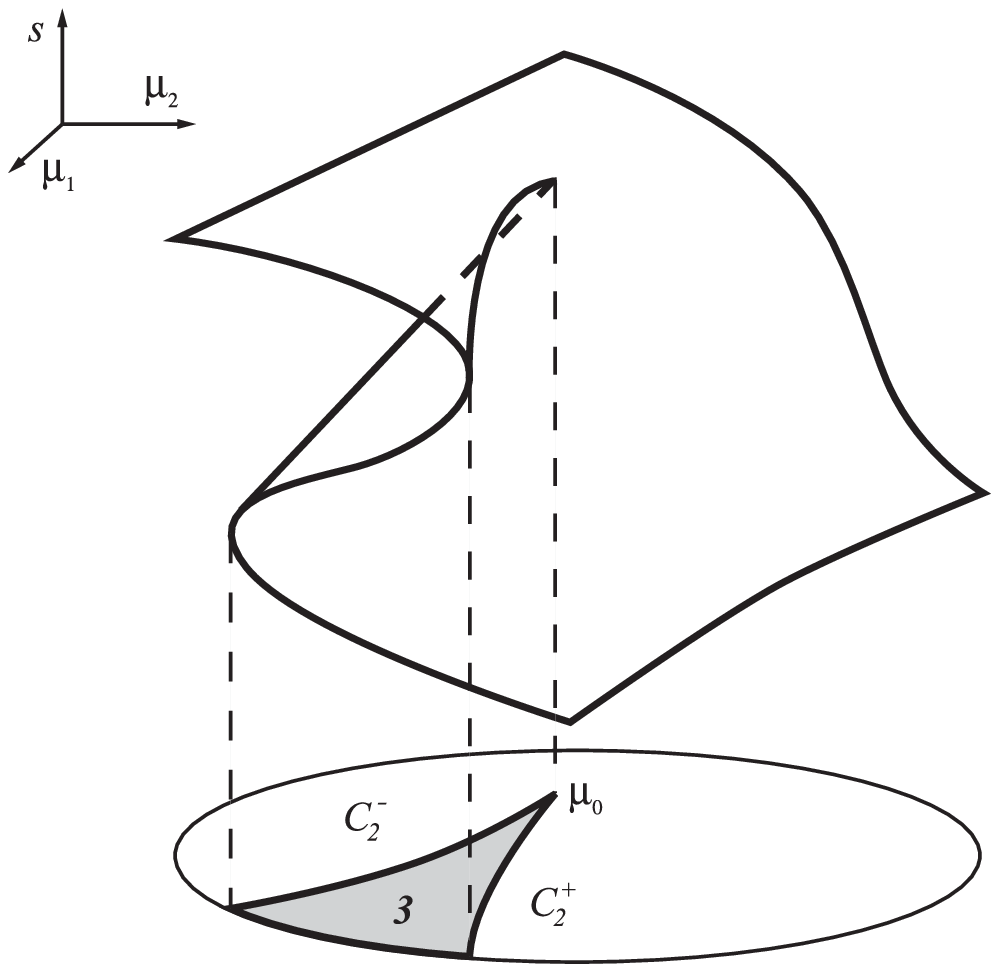}
    \vspace{4mm}
    \par
{\small FIG.~3. The cusp bifurcation surface.}
\end{center}
\end{figure}
    \par
    \medskip
    \textbf{Theorem 3.2.}
    \emph{Suppose that $n\geq3,$ that for $\mbox{\boldmath$\mu$}=
\mbox{\boldmath$\mu$}_0\in\textbf{R}^n$ system (2.2) has
a multiplicity-three limit cycle $L_0,$ that
$d_{\mu_1}(0,\mbox{\boldmath$\mu$}_0)\neq0,$
$d_{r{\mu_1}}(0,\mbox{\boldmath$\mu$}_0)\neq0$
and for $j=2,\ldots,n,$}
    $$
    \Delta_j\equiv\frac{\partial(d, d_r)}
    {\partial(\mu_1,\mu_j)}(0,\mbox{\boldmath$\mu$}_0)\neq0.
    $$
    \par
    \emph{Then given $\varepsilon>0,$ there is a $\delta>0$ and
constants $\omega_j=\pm1$ for $j=2,\ldots,n,$ and there exist
unique functions $h_1(\mu_2,\ldots,\mu_n),$
$h_2(\mu_2,\ldots,\mu_n)$ and $g^{\pm}(\mu_2,\ldots,\mu_n)$
with $h_1(\mu_2^{(0)},\ldots,\mu_n^{(0)})=\mu_1^{(0)},$
$h_2(\mu_2^{(0)},\ldots,\mu_n^{(0)})=\mu_1^{(0)}$ and
$g^{\pm}(\mu_2^{(0)},\ldots,\mu_n^{(0)})=\mu_1^{(0)},$
where $h_1$ and $h_2$ are defined and analytic for
$\vert\mu_j-\mu_j^{(0)}\vert<\delta,$ $j=2,\ldots,n,$
and $g^{\pm}$ are defined and continuous for
$0\leq\sigma_j(\mu_j-\mu_j^{(0)})<\delta$ and analytic
for $0<\omega_j(\mu_j-\mu_j^{(0)})<\delta,$ $j=2,\ldots,n$
such that}
    $$
    C_3: \quad
    \left
    \{
    \begin{array}{rl}
    \mu_{1}=h_{1}(\mu_{2},\ldots,\mu_{n})\\
    \mu_{1}=h_{2}(\mu_{2},\ldots,\mu_{n})\\
    \end{array}
    \right.
    $$
\emph{is an $(n-2)$-dimensional, analytic, cusp bifurcation
surface of multi\-pli\-city-three limit cycles of~(2.2)
through the point $\mbox{\boldmath$\mu$}_0$ and}
    $$
    C_2^{\pm}: \quad
    \mu_{1}=g^{\pm}(\mu_{2},\ldots,\mu_{n})
    $$
\emph{are two $(n-1)$-dimensional, analytic, fold bifurcation
surfaces of multi\-pli\-city-two limit cycles of (2.2) which
intersect in a cusp along
$C_3.$}
    \medskip
    \par
    \textbf{Definition 3.3.}
An analytic surface $C_{4}\subset\textbf{R}^n$ is an
\emph{$(n\!-\!3)$-dimen\-si\-onal swallow-tail bifurcation surface
of multiplicity-four limit cycles of (2.2) through a point
$\mbox{\boldmath$\mu$}_0\in\textbf{R}^n,$} if for all
$\varepsilon>0$ there exists a $\delta>0$ such that for
each $\mbox{\boldmath$\mu$}\in C_{4}$ with
$\|\mbox{\boldmath$\mu$}-\mbox{\boldmath$\mu$}_{0}\|<\delta,$
system (2.2) has a unique multiplicity-four limit cycle
$L_{\mbox{\boldmath$\mu$}}$ in an $\varepsilon$-neighborhood
of $L_{0}$ and system (2.2) undergoes a swallow-tail
bifurcation at $L_{\mbox{\boldmath$\mu$}};$ i.\,e., $C_{4}$ is
the intersection of two $(n\!-\!2)$-dimensional cusp bifurcation
surfaces of multiplicity-three limit cycles $C_{3}^{\pm}$ which
intersect in a cusp along $C_{4};$ furthermore, there are three
$(n\!-\!1)$-dimensional fold bifurcation surfaces of
multiplicity-two limit cycles of (2.2),
$C_{2}^{(i)},$ $i=0,1,2,$ such that $C_{2}^{(0)}$ and
$C_{2}^{(1)}$ intersect in a cusp along $C_{3}^{+},$ $C_{2}^{(0)}$
and $C_{2}^{(2)}$ intersect in a cusp along $C_{3}^{-},$ and
$C_{2}^{(1)}$ and $C_{2}^{(2)}$ intersect along an
$(n\!-\!2)$-dimensional surface on which (2.2) has
two multiplicity-two limit cycles; finally, for
$\|\mbox{\boldmath$\mu$}-\mbox{\boldmath$\mu$}_{0}\|<\delta$
and for $\mbox{\boldmath$\mu$}$ in the swallow-tail region
(shaded in Fig.~4), system (2.2) has four simple limit cycles
in an $\varepsilon$-neighborhood of $L_{0};$ for
$\|\mbox{\boldmath$\mu$}-\mbox{\boldmath$\mu$}_{0}\|<\delta$
and $\mbox{\boldmath$\mu$}$ above the surfaces $C_{2}^{(i)},$
$i=0,1,2,$ system (2.2) has two simple limit cycles in an
$\varepsilon$-neighborhood of $L_{0};$ for
$\|\mbox{\boldmath$\mu$}\!-\!\mbox{\boldmath$\mu$}_{0}\|\!<\!\delta$
and $\mbox{\boldmath$\mu$}$ below the surfaces $C_{2}^{(i)},$
$i=0,1,2,$ system (2.2) has no limit cycles in an
$\varepsilon$-neighborhood of $L_{0}.$ Cf.~Fig.~4.
    \par
\begin{figure}[htb]
\begin{center}
\includegraphics[width=105mm]{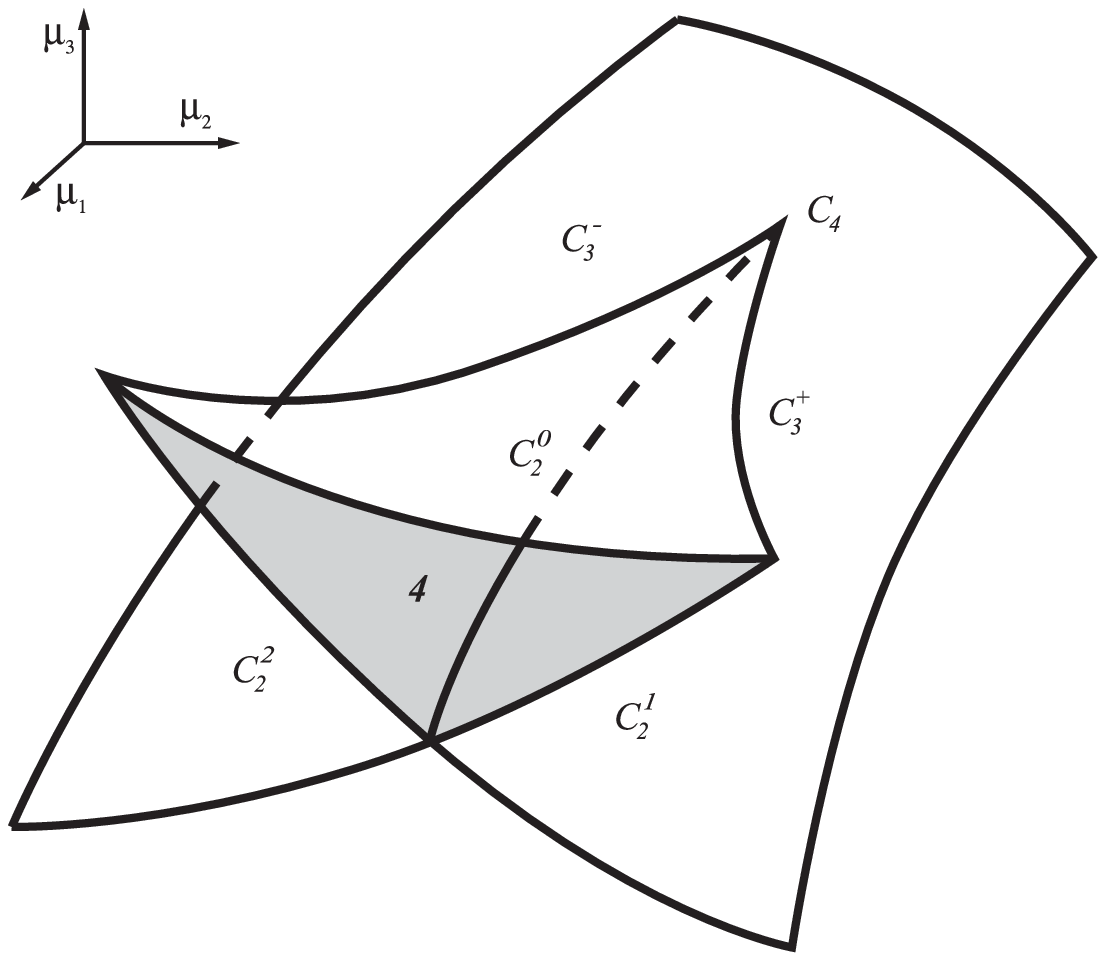}
    \par
{\small FIG.~4. The swallow-tail bifurcation surface.}
\end{center}
\end{figure}
    \par
    \bigskip
    \textbf{Theorem 3.3.}
    \emph{Suppose that $n\geq4,$ that for
$\mbox{\boldmath$\mu$}=\mbox{\boldmath$\mu$}_0\in\textbf{R}^n$
system (2.2) has a multiplicity-four limit cycle $L_0,$
that $d_{\mu_1}(0,\mbox{\boldmath$\mu$}_0)\neq0,$
$d_{r{\mu_1}}(0,\mbox{\boldmath$\mu$}_0)\neq0,$
$d_{rr{\mu_1}}(0,\mbox{\boldmath$\mu$}_0)\neq0,$
and that for $j=2,\ldots,n,$}
    $$
    \frac{\partial(d,d_r)}{\partial(\mu_1,\mu_j)}
    (0,\mbox{\boldmath$\mu$}_0)\neq 0, \quad
    \frac{\partial(d,d_{rr})}{\partial(\mu_1,\mu_j)}
    (0,\mbox{\boldmath$\mu$}_0)\neq 0, \quad
    \frac{\partial(d_r,d_{rr})}{\partial(\mu_1,\mu_j)}
    (0,\mbox{\boldmath$\mu$}_0)\neq0.\\[-2mm]
    $$
    \par
    \emph{Then given $\varepsilon>0,$ there is a $\delta>0$ and
constants $\omega_{jk}=\pm1$ for $j=2,\ldots,n,$ $k=1,2,$ and
there exist unique functions $g_i(\mu_2,\ldots,\mu_n),$
$h_k^{\pm}(\mu_2,\ldots,\mu_n)$ and $F_i(\mu_2,\ldots,\mu_n),$
with $g_i(\mu_2^{(0)},\ldots,\mu_n^{(0)})=
h_k^{\pm}(\mu_2^{(0)},\ldots,\mu_n^{(0)})=
F_i(\mu_2^{(0)},\ldots,\mu_n^{(0)})\!=\!\mu_1^{(0)}\!,$ for
$i=0,1,2$ and $k=1,2,$ where $F_i$ is defined and analytic for
$i=0,1,2,$ and $\vert\mu_j-\mu_j^{(0)}\vert<\delta,$
$j=2,\ldots,n,$ $h_k^{\pm}$ are defined and continuous for
$0\leq\omega_{jk}(\mu_j-\mu_j^{(0)})<\delta$ and analytic for
$0<\omega_{jk}(\mu_j-\mu_j^{(0)})<\delta,$ $j=2,\ldots,n,$
$k=1,2,$ and for $i=0,1,2,$ $g_i$ is defined and analytic in the
cuspidal region between the surfaces
$\mu_1=h_2^{\pm}(\mu_2,\ldots,\mu_n),$ which intersect in a cusp,
and $g_i$ is continuous in the closure of that region, such that}
    $$
    C_4: \quad
    \left
    \{
    \begin{array}{rl}
    \mu_{1}=F_{0}(\mu_{2},\ldots,\mu_{n})\\
    \mu_{1}=F_{1}(\mu_{2},\ldots,\mu_{n})\\
    \mu_{1}=F_{2}(\mu_{2},\ldots,\mu_{n})\\
    \end{array}
    \right.
    $$
\emph{is an $(n-3)$-dimensional, analytic, swallow-tail bifurcation
surface of multi\-pli\-city-four limit cycles of (2.2) through the
point $\mbox{\boldmath$\mu$}_0$ which is the intersection of two
$(n-2)$-dimensional, analytic, cusp bifurcation surfaces of
multiplicity-three limit cycles of (2.2),}
    $$
    C_3^{\pm}: \quad
    \left
    \{
    \begin{array}{rl}
    \mu_{1}=h_{1}^{\pm}(\mu_{2},\ldots,\mu_{n})\\
    \mu_{1}=h_{2}^{\pm}(\mu_{2},\ldots,\mu_{n})\\
    \end{array}
    \right.
    $$
\emph{which intersect in a cusp along $C_4;$ furthermore,
$C_3^{+}=C_2^{(0)} \bigcap C_2^{(1)}$ and $C_3^{-}=C_2^{(0)}
\bigcap C_2^{(2)}$ where for $i=0,1,2,$}
    $$
    C_2^{i}: \quad
    \mu_{1}=g_{i}(\mu_{2},\ldots,\mu_{n})
    $$
\emph{are $(n-1)$-dimensional, analytic, fold bifurcation surfaces
of multiplicity-two limit cycles of (2.2)  which intersect in cusps
along $C_3^{\pm}$ and in an $(n-2)$-dimensional, analytic surface
$C_2^{(1)} \bigcap C_2^{(2)}$ on which (2.2) has two multiplicity-two
limit cycles} (Fig.~4 and Fig.~5).
    \par
\begin{figure}[htb]
\begin{center}
\includegraphics[width=100mm]{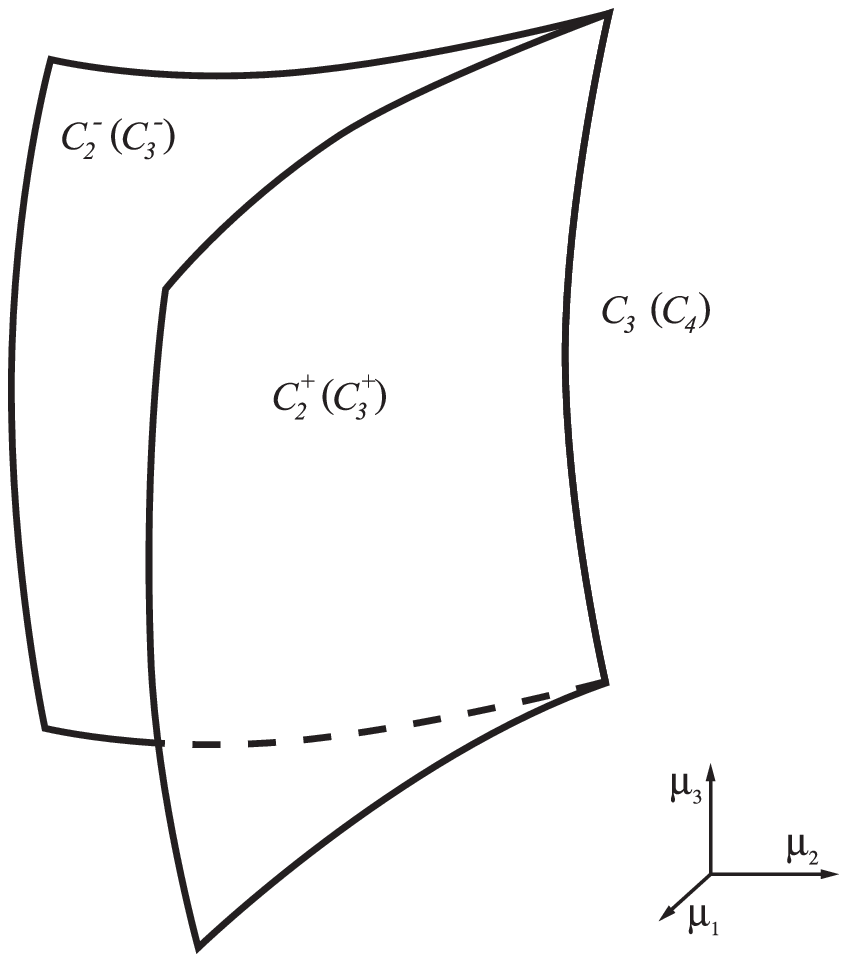}
    \par
{\small FIG.~5. The bifurcation curve (one-parameter family)
of multiple limit cycles.}
\end{center}
\end{figure}
    \medskip
    \par
Based on Theorems~2.3,~2.4, the following ge\-ne\-ra\-li\-za\-tion of
Theorems~3.1\,--\,3.3 can be proved on induction~\cite{Perko}.
    \par
    \medskip
    \textbf{Theorem 3.4.}
    \emph{Given $m\geq2.$ Suppose that $n\geq m,$ that for
$\mbox{\boldmath$\mu$}=\mbox{\boldmath$\mu$}_0\in\textbf{R}^n$
polynomial system~(2.2) has a multiplicity-$m$ limit cycle $L_0,$ that}
    $$
    \frac{\partial d}{\partial\mu_1}\,(0,\mbox{\boldmath$\mu$}_0)\neq0,
        \quad
    \frac{\partial d_r}{\partial\mu_1}\,(0,\mbox{\boldmath$\mu$}_0)\neq0,\;
        \ldots,
        \quad
    \frac{\partial d_r^{(m-2)}}{\partial\mu_1}\,(0,\mbox{\boldmath$\mu$}_0)\neq0,
    $$
\emph{and that}
    $$
    \frac{\partial(d_r^{(i)},d_r^{(j)})}{\partial(\mu_1,\mu_k)}\,
    (0,\mbox{\boldmath$\mu$}_0)\neq0
    \vspace{4mm}
    $$
\emph{for $i,j=0,\ldots,m-2$ with $i\neq j$ and $k=2,\ldots,n.$}
    \par
    \emph{Then given $\varepsilon>0$ there is a $\delta>0$ such that
for $\|\mbox{\boldmath$\mu$}-\mbox{\boldmath$\mu$}_0\|<\delta,$
the system (2.2) has}
    \par
$(1)$ \ \emph{a unique $(n-m+1)$-dimensional analytic surface $C_m$
of multi\-pli\-city-$m$ limit cycles of (2.2) through the point
$\mbox{\boldmath$\mu$}_0;$}
    \par
$(2)$ \ \emph{two $(n-m+2)$-dimensional analytic surfaces $C_{m-1}$
of multiplicity-$(m\!-\!1)$ limit cycles of (2.2) through the point
$\mbox{\boldmath$\mu$}_0$ which intersect in a cusp along $C_m;$}
    \par
    \vspace{-4mm}
    \dots
    \vspace{-2mm}
    \par
$(j)$ \ \emph{exactly $j,$ $(n-m+j)$-dimensional analytic surfaces
$C_{m-j+1}$ of multiplicity-$(m-j+1)$ limit cycles of (2.2) through
the point $\mbox{\boldmath$\mu$}_0$ which intersect pairwise in cusps
along the bifurcation surfaces $C_{m-j+2};$}
    \par
    \vspace{-4mm}
    $\dots$
    \vspace{-2mm}
    \par
$(m-1)$ \ \emph{exactly $(m-1),$ $(n-1)$-dimensional analytic fold
bifurcation surfaces $C_2$ of multiplicity-two limit cycles of (2.2)
through the point $\mbox{\boldmath$\mu$}_0$ which intersect pairwise
in a cusp along the $(n-2)$-dimensional cusp bifurcation surfaces $C_3.$}

    \par
Let us formulate now the Wintner--Perko termination principle~\cite{Perko}
for polynomial system (2.2).
    \par
    \medskip
    \textbf{Theorem 3.5 (Wintner--Perko termination principle).}
    \emph{Any one-para\-me\-ter fa\-mi\-ly of multip\-li\-city-$m$
limit cycles of relatively prime polynomial system (2.2) can be
extended in a unique way to a maximal one-parameter family of
multiplicity-$m$ limit cycles of (2.2) which is either open
or cyclic.}
    \par
\emph{If it is open, then it terminates either as the parameter or
the limit cycles become unbounded; or, the family terminates
either at a singular point of (2.2), which is typically
a fine focus of multiplicity~$m,$ or on a (compound)
separatrix cycle of (2.2) which is also typically of
multiplicity~$m.$}
    \medskip
    \par
The proof of this principle for general polynomial system (2.2)
with a vector parameter $\mbox{\boldmath$\mu$}\in\textbf{R}^n$ parallels
the proof of the pla\-nar termination principle for the system
    $$
    \vspace{1mm}
    \dot{x}=P(x,y,\lambda),
        \quad
    \dot{y}=Q(x,y,\lambda)\\[2mm]
    \eqno(3.1)
    $$
with a single parameter $\lambda\in\textbf{R}$ (see \cite{Gaiko},
\cite{Perko}), since there is no loss of generality in assuming
that system (2.2) is parameterized by a single parameter $\lambda;$
i.\,e., we can assume that there exists an analytic mapping
$\mbox{\boldmath$\mu$}(\lambda)$ of $\textbf{R}$ into $\textbf{R}^n$
such that (2.2) can be written as (3.1) and then we can repeat everything,
what had been done for system (3.1) in~\cite{Perko}. In particular,
if $\lambda$ is a field rotation parameter of (3.1), the following Perko's
theorem on monotonic families of limit cycles is valid; see \cite{Perko}.
    \par
    \medskip
    \textbf{Theorem 3.6.}
    \emph{If $L_{0}$ is a nonsingular multiple limit cycle of (3.1)
for $\lambda=\lambda_{0},$ then $L_{0}$ belongs to a one-parameter family
of limit cycles of (3.1); furthermore:}
    \par
1)~\emph{if the multiplicity of $L_{0}$ is odd, then the family
either expands or contracts mo\-no\-to\-ni\-cal\-ly as $\lambda$
increases through $\lambda_{0};$}
    \par
2)~\emph{if the multiplicity of $L_{0}$ is even, then $L_{0}$
bi\-fur\-cates into a stable and an unstable limit cycle as
$\lambda$ varies from $\lambda_{0}$ in one sense and $L_{0}$
dis\-ap\-pears as $\lambda$ varies from $\lambda_{0}$ in the
opposite sense; i.\,e., there is a fold bifurcation at
$\lambda_{0}.$}

\section{Global bifurcation analysis}

Consider system (1.10). This system has two invariant straight lines: $x=0$
and $y=0.$ Its finite singularities are determined by the algebraic system
   $$
    \begin{array}{l}
x((1-\lambda x)(\alpha x^{2}+\beta x+1)-xy)=0,
    \\[2mm]
y((\delta+\mu y)(\alpha x^{2}+\beta x+1)-x^{2})=0.\\[2mm]
    \end{array}
    \eqno(4.1)
    $$
From (4.1), we have got: two singular points $(0,0)$ and $(0,-\delta/\mu),$
at most two points defined by the condition
    $$
\alpha x^{2}+\beta x+1=0, \quad y=0,
    \eqno(4.2)
    $$
and at most six singularities defined by the system
   $$
    \begin{array}{l}
xy=(1-\lambda x)(\alpha x^{2}+\beta x+1),
    \\[2mm]
y\,(\delta+\mu y)=x\,(1-\lambda x),
    \end{array}
    \eqno(4.3)
    $$
among which we always have the point $(1/\lambda,0).$
    \par
The point $(0,0)$ is always a saddle, but $(1/\lambda,0)$ can be a node or a saddle,
or a saddle-node. The point $(1/\lambda,0)$ can change multiplicity when singular points
enter or exit the first quadrant. In addition, a singular point of multiplicity~2 may appear
in the first quadrant and bifurcate into two singular points. In the case $\beta\geq0$
(respectively, $\beta<0),$ there is a possibility of up to one singular point (respectively,
two singular points) in the open first quadrant \cite{lcr}. If~there exists exactly one simple
singular point in the open first quadrant, then it is an anti-saddle. If~there exists exactly
two simple singular points in the open first quadrant, then the singular point on the left
with respect to the $x$-axis is an anti-saddle and the singular point on the right is a saddle
\cite{lcr}. If a singular point is not in the first quadrant, in consequence, it has no biological
significance.
    \par
To study singular points of (1.10) at infinity, consider the corresponding diffe\-rential equation
    $$
\frac{dy}{dx}=-\frac{y((\delta+\mu y)(\alpha x^{2}+\beta x+1)-x^{2})}
{x((1-\lambda x)(\alpha x^{2}+\beta x+1)-xy)}.
    \eqno(4.4)
    $$
\par
Dividing the numerator and denominator of the right-hand side of (4.4) by $x^{4}$ $(x\neq0)$
and denoting $y/x$ by $u$ (as well as $dy/dx),$ we will get the algebraic equation
    $$
u((\mu/\lambda)u-1)=0, \quad \mbox{where} \quad u=y/x,\\[2mm]
    \eqno(4.5)
    $$
for all infinite singularities of (4.4) except when $x=0$ (the ``ends'' of the $y$-axis),
see \cite{BL}, \cite{Gaiko}. For this special case we can divide the numerator and denominator
of the right-hand side of (4.4) by $y^{4}$ $(y\neq0)$ denoting $x/y$ by $v$ (as well as $dx/dy)$ and consider the algebraic equation
    $$
v^{3}(v-\mu/\lambda)=0, \quad \mbox{where} \quad v=x/y.
    \eqno(4.6)
    $$
The equations (4.5) and (4.6) give three singular points at infinity for (4.4): a~simple node
on the ``ends'' of the $x$-axis, a~triple node on the ``ends'' of the $y$-axis, and a~simple
saddle in the direction of $y/x=\lambda/\mu.$
    \par
To investigate the character and distribution of the singular points in the phase plane,
we have used a method developed in \cite{bg}--\cite{gai10}. The sense of this method is
to obtain the simplest (well-known) system by vanishing some parameters(usually field rotation
parameters) of the original system and then to input these parameters successively one by one
studying the dynamics of the singular points (both finite and infinite) in the phase plane.
    \par
Using the obtained information on singular points and applying a geometric approach developed
in \cite{bg}--\cite{gai10}, we can study the limit cycle bifurcations of system (1.10). This
study will use some results obtained in \cite{lcr}: in particular, the results on the cyclicity
of a singular point of (1.10). However, it is surely not enough to have only these results
to prove the main theorem of this paper concerning the maximum number of limit cycles of
system (1.10).
    \par
Applying the definition of a field rotation parameter \cite{BL}, \cite{Gaiko}, \cite{Perko},
i.\,e., a parameter which rotates the field in one direction, to system (1.10), let us calculate
the corresponding determinants for the parameters $\alpha$ and $\beta,$ respectively:
    $$
\Delta_{\alpha}=PQ'_{\alpha}-QP'_{\alpha}=x^{4}y(y(\delta+\mu y)-x(1-\lambda x)),
    \eqno(4.7)
    $$
    $$
\Delta_{\beta}=PQ'_{\beta}-QP'_{\beta}=x^{4}y(y(\delta+\mu y)-x(1-\lambda x)).\\[4mm]
    \eqno(4.8)
    $$
It follows from (4.7) and (4.8) that on increasing $\alpha$ or $\beta$ the vector field
of (1.10) in the first quadrant is rotated in positive direction (counterclockwise) only
on the outside of the ellipse
    $$
y(\delta+\mu y)-x(1-\lambda x)=0.
    \eqno(4.9)
    $$
Therefore, to study limit cycle bifurcations of system (1.10), it makes sense together
with (1.10) to consider also an auxiliary system (1.11) with a field rotation
parameter~$\gamma\!:$
    $$
\Delta_{\gamma}=P^{2}+Q^{2}\geq0.
    \eqno(4.10)
    $$
    \par
Using system (1.11) and applying Perko's results, we will prove the following theorem.
	\par
    \medskip
\noindent \textbf{Theorem 4.1.}
\emph{System~(1.10) can have at most two limit cycles surrounding one singular point.}
	  \par
    \medskip
\noindent\textbf{Proof.} First let us prove that system (1.10) can have at least two limit cycles.
Begin with system (1.10), where $\alpha=\beta=0$. It is clear that such a cubic system, with two
invariant straight lines, cannot have limit cycles at all \cite{lcr}. Inputting a negative
parameter $\beta$ into this system, the vector field of (1.10) will be rotated in negative
direction (clockwise) at infinity, the structure and the character of stability of infinite
singularities will be changed, and an unstable limit, $\Gamma_{1},$ will appear immediately
from infinity in this case. This cycle will surround a stable anti-saddle (a node or a focus)
$A$ which is in the first quadrant of system (1.10). Inputting a positive parameter $\alpha,$
the vector field of quartic system (1.10) will be rotated in positive direction (counterclockwise)
at infinity, the structure and the character of stability of infinite singularities will be
changed again, and a stable limit, $\Gamma_{2},$ surrounding $\Gamma_{1}$ will appear immediately
from infinity in this case. On further increasing the parameter $\alpha,$ the limit cycles
$\Gamma_{1}$ and $\Gamma_{2}$ combine a semi-stable limit, $\Gamma_{12},$ which then disappears
in a ``trajectory concentration'' \cite{BL}, \cite{Gaiko}. Thus, we have proved that system (1.10)
can have at least two limit cycles; see also \cite{lcr}.
    \par
Let us prove now that this system has at most two limit cycles. The proof is carried out
by contradiction applying Catastrophe Theory; see \cite{Gaiko}, \cite{Perko}. Consider
system (1.11) with three parameters: $\alpha,$ $\beta,$ and $\gamma$ (the parameters
$\delta,$ $\lambda,$ and $\mu$ can be fixed, since they do not generate limit cycles).
Suppose that (1.11) has three limit cycles surrounding the only point $A$ in the first quadrant.
Then we get into some domain of the parameters $\alpha,$ $\beta,$ and $\gamma$ being restricted
by definite con\-di\-tions on three other parameters $\delta,$ $\lambda,$ and $\mu$.
This domain is bounded by two fold bifurcation surfaces forming a cusp bifurcation surface
of multiplicity-three limit cycles in the space of the pa\-ra\-me\-ters $\alpha,$ $\beta,$
and $\gamma$ \cite{Gaiko}, \cite{Perko}.
    \par
The cor\-res\-pon\-ding maximal one-parameter family of multiplicity-three limit cycles
cannot be cyclic, otherwise there will be at least one point cor\-res\-pon\-ding to the
limit cycle of multi\-pli\-ci\-ty four (or even higher) in the parameter space.
Extending the bifurcation curve of multi\-pli\-ci\-ty-four limit cycles through this point
and parameterizing the corresponding maximal one-parameter family of multi\-pli\-ci\-ty-four
limit cycles by the field rotation para\-me\-ter $\gamma,$ according to Theorem~3.6, we will
obtain two monotonic curves of multi\-pli\-ci\-ty-three and one, respectively, which, by the
Wintner--Perko termination principle (Theorem~3.5), terminate either at the point $A$ or on
a separatrix cycle surrounding this point. Since we know at least the cyclicity of the singular
point which is equal to two (see \cite{lcr}), we have got a contradiction with the termination
principle stating that the multiplicity of limit cycles cannot be higher than the
multi\-pli\-ci\-ty (cyclicity) of the singular point in which they terminate.
    \par
If the maximal one-parameter family of multiplicity-three limit cycles is not cyclic,
using the same principle (Theorem~3.5), this again contradicts the cyclicity of $A$
(see \cite{lcr}) not admitting the multiplicity of limit cycles to be higher than two.
This contradiction completes the proof in the case of one singular point in the first
quadrant.
    \par
Suppose that system (1.11) with two finite singularities, a saddle $S$ and an anti-saddle $A,$
has three limit cycles surrounding $A.$ Then we get again into some domain of the parameters
$\alpha,$ $\beta,$ and $\gamma$ bounded by two fold bifurcation surfaces forming a cusp
bifurcation surface of multiplicity-three limit cycles in the space of the pa\-ra\-me\-ters
$\alpha,$ $\beta,$ and $\gamma$ being restricted by definite con\-di\-tions on three other parameters $\delta,$ $\lambda,$ and $\mu$ \cite{Gaiko}, \cite{Perko}.
    \par
The cor\-res\-pon\-ding maximal one-parameter family of multiplicity-three limit cycles
cannot be cyclic, otherwise there will be at least one point cor\-res\-pon\-ding to the
limit cycle of multi\-pli\-ci\-ty four (or even higher) in the parameter space.
Extending the bifurcation curve of multi\-pli\-ci\-ty-four limit cycles through this point
and parameterizing the corresponding maximal one-parameter family of multi\-pli\-ci\-ty-four
limit cycles by the field rotation para\-me\-ter $\gamma,$ according to Theorem~3.6, we will
obtain again two monotonic curves of multi\-pli\-ci\-ty-three and one, respectively, which,
by Theorem~3.5, terminate either at the point $A$ or on a separatrix loop surrounding this
point. Since we know at least the cyclicity of the singular point which is equal to two
(see \cite{lcr}), we have got a contradiction with the termination principle (Theorem~3.5).
    \par
If the maximal one-parameter family of multiplicity-three limit cycles is not cyclic,
using the same principle, this again contradicts the cyclicity of $A$ (see \cite{lcr})
not admitting the multiplicity of limit cycles higher than two. Moreover, it also follows
from the termination principle that a separatrix loop cannot have the multiplicity (cyclicity)
higher than two in this case.
    \par
Thus, we conclude that system~(1.10) cannot have either a multiplicity-three limit cycle or
more than two limit cycles surrounding a singular point which proves the theorem.
\qquad $\Box$

\end{document}